\documentclass{amsart}
\usepackage{enumerate}

\newcommand{\e}{\varepsilon}
\newcommand{\IR}{\mathbb R}
\newcommand{\dens}{\mathrm{dens}}
\newcommand{\w}{\omega}

\newcommand{\Cob}{\mathsf{Cob}}

\newcommand{\pr}{\mathrm{pr}}
\newcommand{\N}{\ensuremath{\mathbb{N}}}
\newcommand{\clo}{\overline}

\newtheorem{theorem}{Theorem}[section]
\newtheorem{lemma}[theorem]{Lemma}
\newtheorem{corollary}[theorem]{Corollary}
\newtheorem{definition}[theorem]{Definition}
\newtheorem{construction}[theorem]{Construction}
\newtheorem{proposition}[theorem]{Proposition}
\newtheorem{problem}[theorem]{Problem}

\newtheorem{example}[theorem]{Example}

\title{Nonseparably connected complete metric spaces}
\author{T.~Banakh, M.~Vovk, M.~R.~W\'ojcik}
\address{Taras Banakh: Ivan Franko Lviv National University, Lviv, Ukraine, and\newline
Unwersytet Humanistyczno-Przyrodniczy im. Jana Kochanowskiego, Kielce, Poland}
\email{tbanakh@yahoo.com}
\address{M.~Vovk: National University ``Lvivska Politechnika", Lviv, Ukraine}
\address{Micha{\l} Ryszard W\'ojcik:\newline Department of Mathematics, University of Louisville, Louisville, USA}
\email{michal.ryszard.wojcik@gmail.com}

\begin{document}

\maketitle


\begin{abstract}
A topological space is {\em nonseparably connected} if it is connected
but all of its connected separable subspaces are singletons. 
We show that each connected first countable space is the image of a nonseparably connected complete metric space
under a continuous monotone hereditarily quotient map.
\end{abstract}


\section{Introduction}
In this paper we solve the problem of constructing a nonseparably connected complete metric space, posed in \cite[Problem 2]{MW2}.
A topological space is {\em separably connected} if any two of its points lie in a connected separable subspace.
On the other hand, a topological space is {\em nonseparably connected} if it is connected
but all of its connected separable subspaces are singletons. 

The first example of a nonseparably connected metric space was constructed by Pol
in 1975, \cite{Pol}. Another example was given by Simon in 2001, \cite{Sim}.
In 2008, Morayne and W\'ojcik obtained a nonseparably
connected metric group as a graph of an additive function from the real line,
\cite{MW2}. None of these nonseparably connected spaces are completely metrizable.

\begin{definition}{\em
The {\em separablewise component} of a point $x_0$ of a topological space $X$
is the union of all separable connected subsets containing $x_0$.
A space is {\em separably connected} if it has only one separablewise component.
A space is {\em nonseparably connected} if it is connected but all of its
separablewise components are singletons.
}\end{definition}

Let $A$ be a separablewise component of a point $x$ in a sequential space $X$.
We will show that $A$ is closed.
If $a\in\clo A$, then there is a sequence of points $(a_n)_{n=1}^{\infty}\subset A$
converging to $a$.
There is a sequence of connected separable sets $A_n$ with $x\in A_n$ and $a_n\in A_n$.
Notice that the set $\{a\}\cup\bigcup\{A_n\colon n\in\N\}$ is connected and separable.
So it must be contained in $A$. In particular, $a\in A$.
\\\\
This is a long paper so we decided to precede it with a short overview
to acquaint the reader with our main results and
to help navigate through the different sections.

For any connected metric space $(X,d)$ we construct a connected complete metric space $\Cob(X,d)$
called the cobweb over $(X,d)$ and define a continuous monotone hereditarily quotient
compression map $\pi\colon\Cob(X,d)\to X$ whose fibers coincide
with the separablewise components of $\Cob(X,d)$.

The cobweb space $\Cob(X,d)$ is embedded in $\Gamma(X)$, called the complete oriented graph over $X$,
which is simply a graph in which any two distinct vertices $x,u$ are joined by two separate edges
denoted $[x,u]$ and $[u,x]$ in such a way that the edges carry the Euclidean metric
and become arcs and the distance between points on different edges is measured
as the shortest distance of traveling along the edges.

Next we iterate the cobweb construction
taking the cobweb over some connected metric space $(X,d)$ and then the cobweb over this cobweb and so on.
We take the inverse limit of this sequence of cobwebs, $\Cob^\omega(X,d)$,
to kill off all separable connected subsets
and end up with a nonseparably connected complete metric space.

Finally, we notice that the cobweb construction works for all distance spaces $(X,d)$,
which is a huge generalization of metric spaces, and study such spaces in some detail
to formulate our strongest result.

As an application of the cobweb construction we present a non-constant continuous real-valued function
defined on a connected complete metric space that has a local minimum or a local maximum at every point.

Since the notion of a quotient map plays a crucial role in our constructions
we decided to collect all the necessary definitions and basic theorems into an appendix.

\section{The cobweb over a distance space}\label{cobweb}

Our fundamental tool is the cobweb construction which was originally conceived for metric spaces,
but turns out to work for a much broader class of spaces.

\begin{definition}{\em
For any function $d\colon X\times X\to[0,\infty)$
such that $d(x,x)=0$ for all $x\in X$,
we call the ordered pair $(X,d)$ a {\em distance space}.

Let each set of the form $B_d(x,r)=\{z\in X\colon d(x,z)<r\}$
be called a {\em ball of radius $r>0$ centered at $x\in X$}.

We equip each distance space $(X,d)$ with the {\em topology generated by the distance function} $d$
so that a set $E\subset X$ is defined to be {\em open in }$(X,d)$
if for every point $x\in E$ there is an $r>0$ such that $B_d(x,r)\subset E.$

A distance space $(X,d)$ is called {\em well-behaved} if
$x\in Int(B_d(x,r))$ for all $x\in X$ and $r>0$.
}\end{definition}

It is enough to note that every metric space is a well-behaved distance space
to go through our basic constructions.
Later we will characterize distance spaces and well-behaved distance spaces
when formulating the strongest result in Section \ref{Eco}.

\begin{definition}
Let $\kappa$ be a cardinal number.
A {\em hedgehog with $\kappa$ spikes, each of length $\e>0$},
is the space $H=\{(0,0)\}\cup(\kappa\times(0,\e])$
equipped with the metric $\rho$ given by
\[\rho((x,t),(u,s))=\begin{cases}
|t-s|&\text{if}\ \ x=u,\\
\ t+s&\text{if}\ \ x\not=u.
\end{cases}\]
\end{definition}

It is easy to see that a hedgehog is an arcwise connected complete metric space.

\begin{definition}{\em
{\em The complete oriented graph $\Gamma(X)$
over a set $X$ of vertices}
is defined to be the union
$$\Gamma(X)=\bigcup\{[x,y]\colon x,y\in X,x\not=y\}$$
of all oriented edges $[x,y]$ with distinct vertices $x,y\in X$,
so that $$[x,y]\cap[y,x]=\{x,y\},$$
and the oriented edge $[x,y]$ is defined as
$$[x,y]=\{x,y\}\cup\{(x,y,t)\colon t\in(0,1)\}.$$
}\end{definition}

\begin{construction}{\em
We are going to define a complete metric $\rho$ on $\Gamma(X)$ such that
\begin{enumerate}
\item each edge $[x,y]$ is isometric to the unit interval $[0,1]$,
\item each edge without end points $[x,y]\setminus\{x,y\}$ is an open set,
\item if two points $a,b$ lie on disjoint edges then $\rho(a,b)>1$.
\end{enumerate}
}\end{construction}
\begin{proof}
First, we equip each edge with the euclidean metric
by defining an auxiliary function $r$ on a subset of $\Gamma(X)\times\Gamma(X)$.
For all $x,y\in X,x\not=y$ and all $t,s\in(0,1)$,
let $$r(x,x)=0,\ r(x,y)=1,$$
$$r((x,y,t),(x,y,s))=|t-s|,$$
$$r(x,(x,y,t))=r((x,y,t),x)=t,$$
$$r(y,(x,y,t))=r((x,y,t),y)=1-t.$$
Next, we extend this auxiliary function to the whole $\Gamma(X)\times\Gamma(X)$.
For any points $a,b\in\Gamma(X)$,
let $\rho(a,b)$ be the infimum over finite sums
$$\sum_{i=1}^n r(a_i,a_{i-1})$$
where $\{a_0,\ldots,a_n\}\subset X$ with $a=a_0$, $a_n=b$ and
consecutive points $a_i,a_{i-1}$ belong to the same edge:
$(\forall i)(\exists x,y\in X)\ \{a_i,a_{i-1}\}\subset[x,y].$

It is clear that $\rho$ is a metric on $\Gamma(X)$ from the way it is defined.

It is easy to see that each edge $[x,y]$ is isometric to $[0,1]$.

The set $[x,y]\setminus\{x,y\}$ is open because
it is equal to the open ball of radius $1/2$
centered in the middle of $[x,y]$, that is at $(x,y,1/2)$.

The distance between any two points on disjoint edges is greater than one
because to travel from one such point to another along the edges
we must move along at least one whole edge.

Each closed ball of radius $1/2$ centered at a vertex $x\in X$
is isometric to a hedgehog with $2(|X|-1)$ spikes of length $1/2$.
In particular, such closed balls are complete metric spaces.

The complete oriented graph $(\Gamma(X),\rho)$ is a complete metric space.
Indeed, let $x_n$ be a Cauchy sequence such that
$\rho(x_n,x_m)<1/4\text{ for all }n,m\in\mathbb N.$
If $\rho(x_n,u)\ge 1/4$ for all $n\in\mathbb N$ and all $u\in X$,
then our sequence is contained in one of the edges.
Otherwise, there is an index $k$ and a vertex $u\in X$ such that
$\rho(x_k,u)<1/4$. Then $\rho(x_n,u)\le\rho(x_n,x_k)+\rho(x_k,u)\le 1/2$ for all $n\in\mathbb N$.
This means that our sequence
is contained in the closed ball of radius $1/2$ centered at $u$,
which is a complete metric space.
In both cases, our Cauchy sequence converges.
\end{proof}

\begin{definition}{\em
For a distance space $(X,d)$,
we construct a metric space $\Cob(X,d)$,
the so called {\em cobweb over the distance space $(X,d)$},
as a special subspace of $(\Gamma(X),\rho)$ and the so called {\em compression map}
$\pi\colon\Cob(X,d)\to X$ in the following way:

Let $d_1=\min\{d,1/2\}$.
For distinct vertices $a,b\in X$, let $a_b\in[a,b]$
be the unique point satisfying $$\rho(a_b,b)=d_1(b,a).$$
Note that $a_b\not=b\iff d(b,a)>0$.
Let $[a,a_b]$ be the subarc of $[a,b]$ given by
$$[a,a_b]=\{z\in[a,b]\colon\rho(a,z)\le 1-d_1(b,a)\}.$$
Let $$Z_a=\bigcup\big\{[a,a_b]\setminus\{b\}\ \colon\ b\in X\setminus\{a\}\big\}.$$
Let $$\Cob(X,d)=\bigcup_{a\in X}Z_a.$$

Let the so called {\em compression map}
$\pi\colon\Cob(X,d)\to X$ be given by $$\pi(Z_a)=\{a\}\text{ for all }a\in X.$$
}\end{definition}

The next theorem describes the basic properties of the cobweb over a distance space
without any reference to the topology of the distance space.
In Section \ref{applications} devoted to the applications of the cobweb construction,
we will refer to these properties only.

\begin{theorem}\label{summary}
Let $(X,d)$ be a distance space.
Let $\pi\colon\Cob(X,d)\to X$ be the compression map.
Then for all $x\in X$ and all $r\in(0,\frac{1}{2})$,
\begin{enumerate}
\item $X\subset\Cob(X,d)$,
\item $\pi(x)=x$,
\item $\pi^{-1}(x)$ is arcwise connected,
\item $\pi^{-1}(x)\setminus\{x\}$ is open in $\Cob(X,d)$,
($\pi$ is locally constant on $\Cob(X,d)\setminus X$),
\item $(\Cob(X,d),\rho)$ is a complete metric space,
\item $\pi(B_\rho(x,r))=B_d(x,r)$,
\item $|X|\le\dens(\Cob(X,d))\le|\Cob(X,d)|=\mathfrak c|X|$.
\end{enumerate}
\end{theorem}
\begin{proof}
(1)--(4) Evident from the definition.

(5) Notice that $\Cob(X,d)$ is a closed subset of $(\Gamma(X),\rho)$
because it is obtained by taking away selected open intervals from some of the edges
so that in effect it is the intersection of a family of closed subsets.
Since $(\Gamma(X),\rho)$ is complete,
$(\Cob(X,d),\rho)$ is a complete metric space.

(6) If $x\in B_d(a,r)$, then $\rho(a,x_a)=d(a,x)<r$.
So $x_a\in B_\rho(a,r)$ and $x=\pi(x_a)$.
Thus $x\in \pi(B_\rho(a,r))$.

On the other hand,
if $x\in \pi(B_\rho(a,r))$, we have a point $z\in\Cob(X,d)$
such that $\rho(a,z)<r$ and $\pi(z)=x$.
Consequently, $z\in[x,x_u]$ for some $u\in X\setminus\{x\}$.
We may assume that $x\not=a$ because otherwise there is nothing to prove.
Now, supposing that $u\not=a$, we have three distinct vertices $a,x,u$
and a point $z$ lying on the edge $[x,u]$,
which means that $\rho(a,z)\ge 1$.
This contradiction shows that $u=a$ and thus $z\in[x,x_a]$.
Now, the point $x_a$ lies between $z$ and $a$ on the edge $[x,a]$,
so that $d(a,x)=\rho(x_a,a)\le\rho(a,z)<r$,
showing that $x\in B_d(a,r)$.

(7) Notice that each fiber of the compression map $\pi$ contains a hedgehog
with $|X|-1$ spikes. Thus $|\Cob(X,d)|=\mathfrak c|X|$.
The nonempty set $\pi^{-1}(x)\setminus\{x\}$ is open in $\Cob(X,d)$ for each $x\in X$,
so $|X|\le\dens(\Cob(X,d))$.
\end{proof}

We encourage the reader to check that
if $(X,d)$ is a distance space such that
$d(a,b)>0$ for all distinct $a,b\in X$,
then
every fiber of the compression map is homeomorphic
to a hedgehog with $|X|-1$ spikes;
and the set of points
at which the space $\Cob(X,d)$ is locally connected
is equal to $\Cob(X,d)\setminus X$.

Refer to Section \ref{appendix}
for definitions of some terms used in the following theorem.

\begin{theorem}\label{summary2}
For any distance space $(X,d)$,
the compression map $\pi\colon\Cob(X,d)\to X$
is a continuous monotone quotient surjection
and $\Cob(X,d)$ is connected iff $(X,d)$ is connected.
Moreover, $\pi$ is hereditarily quotient iff $(X,d)$ is well-behaved.
\end{theorem}
\begin{proof}
Since $\pi(x)=x$ for every $x\in X$, $\pi$ is a surjection.

Since the set $\pi^{-1}(x)\setminus\{x\}$ is open for every $x\in X$,
the map $\pi$ is locally constant at each point $z\in\Cob(X,d)\setminus X$.
In particular, it is continuous at these points.
That $\pi$ is continuous at each $x\in X$ follows directly from
the inclusion $$\pi(B_\rho(x,r))\subset B_d(x,r)$$
and the definition of the topology on $X$.
So $\pi$ is continuous.

To show that $\pi$ is quotient, take any $x\in A$ such that 
$\pi^{-1}(A)\text{ is open in }\Cob(X,d)$. We need to show that $A\text{ is open in }X.$
Since $\pi(x)=x$, we have $x\in \pi^{-1}(A)$.
Since $\pi^{-1}(A)$ is open, there is an $r>0$ such that $B_\rho(x,r)\subset \pi^{-1}(A)$.
Then $$B_d(x,r)\subset \pi(B_\rho(x,r))\subset A$$
showing that $A$ is open in $(X,d)$.

The compression map is monotone because each fiber
$\pi^{-1}(x)$ is arcwise connected.
We have just proved that $\pi$ is quotient.
So by Lemma \ref{ZconnectedifXconnected},
if $(X,d)$ is connected, then $\Cob(X,d)$ is connected.
On the other hand, if $\Cob(X,d)$ is connected,
then $(X,d)$ is connected because the compression map is continuous.

Suppose now that $(X,d)$ is well-behaved.
To show that $\pi$ is hereditarily quotient,
take any $x\in X$ and any open set $U\subset\Cob(X,d)$ such that $\pi^{-1}(x)\subset U$.
Then $x\in\pi^{-1}(x)\subset U$,
and since $U$ is open, $x\in B_\rho(x,r)\subset U$ for some $r>0$.
Finally, $x\in Int(B_d(x,r))\subset B_d(x,r)\subset\pi(B_\rho(x,r))\subset\pi(U)$
and thus $x\in Int(\pi(U))$.

Suppose now that the compression map $\pi\colon\Cob(X,d)\to X$
is hereditarily quotient.
To show that $(X,d)$ is well-behaved,
take any $x\in X$ and $r>0$.
Let $$U=B_\rho(x,r)\cup\big(\pi^{-1}(x)\setminus\{x\}\big).$$
Notice that $U$ is open in $\Cob(X,d)$ and that $\pi^{-1}(x)\subset U$.
Since $\pi$ is hereditarily quotient, $x\in Int(\pi(U))$.
But $\pi(U)=B_d(x,r)$, so $x\in Int(B_d(x,r))$.
\end{proof}

To describe the separablewise components of the cobweb space now
and to obtain economical metrics later,
we will make use of the fact that
the compression map is locally constant except on a metrically discrete subset.
Naturally,
the cardinality of the image $f(X)$ of a locally constant function $f\colon X\to Y$
does not exceed the density of the domain, $|f(X)|\le\dens(X)$,
and the cardinality of a metrically discrete space does not exceed its density.

\begin{theorem}\label{sepcomp}
Let $(X,d)$ be a distance space.
Let $\pi\colon\Cob(X,d)\to X$ be the compression map.
Then $|\pi(A)|\le\dens(A)$ for any $A\subset\Cob(X,d)$.
Consequently, for any metric space $(X,d)$,
the fibers of $\pi$ coincide with the separablewise components
of $\Cob(X,d)$ which in turn coincide with the arcwise components.
\end{theorem}
\begin{proof}
Let $A\subset\Cob(X,d)$.
Notice that the restriction $\pi|(A\setminus{X})$ is locally constant,
so 
$|\pi(A)|\le|\pi(A\cap X)|+|\pi(A\setminus{X})|\le|A\cap X|+\dens(A\setminus{X}).$
On the other hand, $|A\cap X|\le\dens(A\cap X)$ because $A\cap X$ is metrically discrete,
$\rho(a,b)=1$ for distinct $a,b\in A\cap X$.
Thus $|\pi(A)|\le\dens(A\cap X)+\dens(A\setminus {X})=\dens(A).$

Let $E$ be a connected separable subset of $\Cob(X,d)$.
Since $\pi$ is continuous, $\pi(E)$ is connected.
Since $|\pi(E)|\le\dens(E)=\aleph_0$, $\pi(E)$ is countable.
Now, if the connected countable set $\pi(E)$ is contained
in the metric space $X$, it must be a singleton,
which means that $E$ lies in one of the fibers of $\pi$.
Recall that the fibers are arcwise connected.
\end{proof}


The following theorem reveals the surprising abundance
of connected complete metric spaces that fail to be separably connected.

\begin{theorem}\label{abundance}
For every connected metric space $(X,d)$
the cobweb $\Cob(X,d)$ is
a connected complete metric space whose separablewise components
form a quotient space homeomorphic to $X$.
\end{theorem}
\begin{proof}
By Theorem \ref{summary2}, $\Cob(X,d)$ is a connected complete metric space,
whose separablewise components, by Theorem \ref{sepcomp},
coincide with the fibers of the compression map,
which in turn form a quotient space homeomorphic to $X$,
as may be easily verified by using the following distance function
$$D\big(\pi^{-1}(a),\pi^{-1}(b)\big)=d(a,b)=d(b,a)$$
to generate the topology of the quotient space
formed by the fibers.
(In fact, if $(X,d)$ is any symmetric distance space
then the topology of the quotient space formed by the fibers of
the compression map is generated by
this distance function.)
\end{proof}

Notice that
in this way we encode the metric space $(X,d)$ as a completely different object
$\Cob(X,d)$ from which we can extract the original space by a topological operation.

Let us make two final remarks to conclude this section.

Recall that a connected, locally connected complete metric space
must be arcwise connected, \cite[6.3.11]{En}.
The cobweb over a connected metric space
is not locally connected, although it is locally connected except on
a metrically discrete subset.
This illustrates how important it is
to assume that the space is locally connected at each point
if we want to conclude that it is arcwise connected.

For any distance space $(X,d)$
the separablewise component of a point $x\in X\subset\Cob(X,d)$
coincides with the fiber $\pi^{-1}(x)$ if and only if
each countable connected subset $C$ of the distance space $X$ with $x\in C$
is equal to the singleton $\{x\}$.
This brings us to the following notion.
We say that a topological space $X$ is {\em functionally Hausdorff}
if for any two distinct points $a,b\in X$ there is a continuous
function $f\colon X\to\mathbb R$ with $f(a)\not=f(b)$.
Every metric space $(X,d)$ is functionally Hausforff,
because for any distinct points $a,b\in X$, the function
$f\colon X\to[0,1]$ given by $f(x)={d(x,a)}/({d(x,a)+d(x,b)})$
is continuous and $f(a)=0$, $f(b)=1$.
Each connected subset $E$ of a functionally Hausforff space is either a singleton or contains
a set of cardinality $\mathfrak c$.
Indeed, let $a,b\in E$ with $f(a)<f(b)$.
Then the connected set $f(E)$ must contain the interval $[f(a),f(b)]$.
In conclusion, Theorem \ref{sepcomp}
could be stated for any functionally Hausdorff distance space $(X,d)$
and not just for metric spaces.
(The lexicographical square with the order topology
is a connected compact Hausdorff first countable space.
Being normal, it is functionally Hausdorff.
Being first countable, it may serve as an example of
a connected functionally Hausdorff well-behaved distance space
that is not metrizable.)

\section{The iterated cobweb functor}\label{iterated}

\begin{definition}\label{cobomega}{\em
Given a distance space $(X,d)$,
we define by induction a sequence of iterated cobweb spaces over $(X,d)$:
let $\Cob^1(X,d)=\Cob(X,d)$ be equipped with the natural cobweb metric $\rho_1$
induced from the complete graph over $X$,
and let $\Cob^{n+1}(X,d)=\Cob(\Cob^n(X,d),\rho_n)$
be equipped with the natural cobweb metric $\rho_{n+1}$
induced from the complete graph over $\Cob^n(X,d)$.

Let $\pi_n^{n+1}\colon\Cob^{n+1}(X,d)\to\Cob^n(X,d)$
denote the appropriate compression maps.

Let
$$\Cob^\omega(X,d)=\big\{(x_n)_{n=1}^\infty\in\prod_{n\in\N}\Cob^n(X,d)
\ \colon\ \pi_n^{n+1}(x_{n+1})=x_n\text{ for all }n\in\N\big\}$$
be the inverse limit of the sequence of iterated cobweb spaces over $(X,d)$.
}\end{definition}

\begin{theorem}\label{summary10}
For any distance space $(X,d)$,
$\Cob^\omega(X,d)$ is a completely metrizable space
that is connected iff $(X,d)$ is connected,
and its separablewise components are singletons.
\end{theorem}
\begin{proof}
$\Cob^\omega(X,d)$ is completely metrizable
because all the factor spaces $\Cob^n(X,d)$ are completely metrizable.

Suppose that $(X,d)$ is connected.
Then by Theorem \ref{summary2},
the spaces $\Cob^n(X,d)$ are connected and
the compression maps
$\pi_n^{n+1}\colon\Cob^{n+1}(X,d)\to\Cob^n(X,d)$
are continuous monotone hereditarily quotient surjections for all $n\in\N$.
Therefore, by Puzio's Theorem \ref{Puzio}, $\Cob^\omega(X,d)$ is connected.

Let us write $\pr_n((x_k)_{k=1}^\infty)=x_n$ for $n\in\N$.

Let $A\subset\Cob^\omega(X,d)$ be connected and separable.
Fix $n\in\N$.
The projection $\pr_{n+1}(A)$ is connected and separable
so, by Theorem \ref{sepcomp}, it lies in one of the fibers of
$\pi_n^{n+1}$.
Thus $\pi^{n+1}_n(\pr_{n+1}(A))$ is a singleton.
But $\pi^{n+1}_n(\pr_{n+1}(A))=\pr_n(A)$.
So $\pr_n(A)$ is a singleton for each $n\in\N$,
which means that $A$ is a singleton.
\end{proof}

\begin{corollary}
For every connected distance space $(X,d)$
having at least two points,
$\Cob^\omega(X,d)$ is a nonseparably connected
completely metrizable space.
\end{corollary}

Recall that a topological space is {\em punctiform}
if all of its connected compact subsets are singletons.
For example, any Bernstein subset of the euclidean plane is
a connected punctiform metric space.
A separable connected punctiform complete metric space was constructed
by Kuratowski and Sierpi{\'n}ski in 1922, \cite{connexes_punctiformes}.
A variation of their idea was presented in \cite{WPhD}.
Our nonseparably connected complete metric spaces
are new examples of connected punctiform complete metric spaces.

Using the notion of an economical metric we show that all separable subsets of our nonseparably connected spaces
are in fact zero-dimensional.

\begin{definition}\label{econo}
{\em
Given a metric space $(X,d)$, we say that the metric $d$ is {\em economical} if
$\mathrm{card}(\{d(a,b)\colon a,b\in A\})=|d(A\times A)|\le \dens(A)=\min\{|D|\colon D\subset A\subset\clo D\}$
for any infinite subset $A\subset X$.
We say that a topological space is {\em economically metrizable}
if its topology can be generated by an economical metric.
}\end{definition}

\begin{proposition}\label{eco-p1}
If $X$ is an economically metrizable space,
then each subspace $A\subset X$ of density $\dens(A)<\mathfrak c$ is zero-dimensional.
Consequently, each connected economically metrizable space is nonseparably connected.
\end{proposition}
\begin{proof}
Take any $a\in A$ and any $r>0$.
If $(\forall\e\in(0,r))(\exists b\in A)(d(a,b)=\e)$,
then $|d(\{a\}\times A)|\ge|(0,r)|=\mathfrak c$.
This contradiction shows that there is an $\e\in(0,r)$
such that $d(a,b)\not=\e$ for all $b\in A$.
So the open ball of radius $\e$ centered at $a$ $B(a,\e)\cap A$
is relatively clopen in $A$ and contained in $B(a,r)$.
This means that $A$ is zero-dimensional.
\end{proof}

\begin{theorem}\label{complete_economical}
Let $(X,d)$ be a distance space.
Then $$\rho_\infty\big((x_n)_{n=1}^\infty,(u_n)_{n=1}^\infty\big)=
\max\Big\{\cfrac{\rho_n(x_n,u_n)}{n}\colon{n\in\N}\Big\}$$
is a complete economical metric for $\Cob^\omega(X,d)$.
Consequently, for every connected distance space $(X,d)$,
$(\Cob^\omega(X,d),\rho_\infty)$ is a connected
complete economical metric space and its separable subsets
are zero-dimensional.
\end{theorem}
\begin{proof}
Notice that in the definition above we have $\rho_n\le 2$ for all $n\in\N$
because the complete graph metrics are always bounded by two.

Notice that the distance function $\rho_\infty$
is a product metric on $\prod_{n\in\N}\Cob^n(X,d)$.
Since all the factor metric spaces $(\Cob^n(X,d),\rho_n)$ are complete,
so is the product metric space $(\prod_{n\in\N}\Cob^n(X,d),\rho_\infty)$.
Now, $\Cob^\omega(X,d)$, as an inverse limit, is a closed subset and thus
a complete metric space.

Let us write $\pr_n((x_k)_{k=1}^\infty)=x_n$ for $n\in\N$.

Notice that for any $a,b\in\prod_{n\in\N}\Cob^n(X,d)$ we have
$$\rho_\infty(a,b)\in\bigcup\{n^{-1}\rho_n(\pr_n(a),\pr_n(b))\colon n\in\N\}.$$
Therefore, if $A\subset\prod_{n\in\N}\Cob^n(X,d)$, we have
$$|\rho_\infty(A\times A)|\le\sum_{n\in\N}|\rho_n(\pr_n(A)\times\pr_n(A))|\le
\sum_{n\in\N}|\pr_n(A)|^2.$$

To show that $(\Cob^\omega(X,d),\rho_\infty)$ is economical,
take any infinite $A\subset\Cob^\omega(X,d)$.
By Theorem \ref{sepcomp},
$|\pi_n^{n+1}(E)|\le\dens(E)$
for all $n\in\N$ and $E\subset\Cob^{n+1}(X,d)$.

Therefore, for all $n\in\N$,
$$|\pr_n(A)|=|\pi_n^{n+1}(\pr_{n+1}(A))|\le\dens(\pr_{n+1}(A))\le\dens(A).$$
Thus, since $\dens(A)$ is infinite,
$$|\rho_\infty(A\times A)|\le\sum_{n\in\N}|\pr_n(A)|^2
\le\sum_{n\in\N}(\dens(A))^2=\dens(A).$$
\end{proof}

\section{More on distance spaces}\label{distance}

Obviously, every metric space is a well-behaved distance space.
Moreover, the topology of every first countable space
can be generated by a well-behaved distance function
and conversely every well-behaved distance space is first countable.
Furthermore, each distance space (not necessarily well-behaved)
is weakly first countable in the sense of Arkhangelskii \cite{Arh},
and conversely the topology of every weakly first countable space
can be generated by a distance function.

\begin{definition}[Arkhangelskii \cite{Arh}]\label{weakly}
A topological space $X$ is {\em weakly first countable}
if to each point $x\in X$ we can assign a decreasing sequence
$(B_n(x))_{n\in\w}$ of subsets of $X$ that contain $x$ so that a subset
$U\subset X$ is open if and only if for each $x\in U$ there is $n\in\w$ with $B_n(x)\subset U$.
\end{definition}

\begin{lemma}\label{beobachtung}
Let $X$ be an arbitrary set.
Let $E_n(x)\subset X$, $n\in\N,x\in X$, be arbitrary sets
such that $x\in E_{n+1}(x)\subset E_n(x)$ for all $x\in X$ and $n\in\N$.
Let $d\colon X\times X\to[0,1]$ be given by
$$d(x,y)=\inf\{1/n\colon y\in E_n(x)\}.$$
Then $$E_{n+1}(x)\subset B_d(x,1/n)\subset E_n(x)$$
for all $x\in X$ and $n\in\N$.
\end{lemma}
\begin{proof}
Fix $x\in X$ and $n\in\mathbb N$.

If $y\in B_d(x,1/n)$, then $d(x,y)<1/n$ and so there is a $k\in\mathbb N$
such that $1/k<1/n$ and $y\in E_k(x)$.
Since $k>n$, $y\in E_k(x)\subset E_n(x)$.
In effect, $B_d(x,1/n)\subset E_n(x).$

If $y\not\in B_d(x,1/n)$, then $d(x,y)\ge 1/n$, and thus
$(\forall k\in\mathbb N)(y\in E_k(x)\implies 1/k\ge 1/n).$
So $y\not\in E_{n+1}(x)$.
In effect, $X\setminus B_d(x,1/n)\subset X\setminus E_{n+1}(x)$.
\end{proof}

\begin{theorem}\label{1stcountablewell}
A topological space is weakly first countable if and only if
its topology is generated by a distance function.
A topological space is first countable if and only if
its topology is generated by a well-behaved distance function.
\end{theorem}
\begin{proof}
Notice that for a distance space $(X,d)$, the family of sets
$$\{B_d(x,1/n)\colon x\in X,n\in\mathbb N\}$$
witnesses that it is weakly first countable.
If the distance function $d$ is well-behaved, then
the family of open sets
$$\{Int(B_d(x,1/n))\colon x\in X,n\in\mathbb N\}$$
witnesses that it is first countable.

Let $X$ be a weakly first countable space. We have sets
$\{B_n(x)\colon x\in X,n\in\mathbb N\}$
as in Definition \ref{weakly}.
Let $d\colon X\times X\to[0,1]$ be given by
$$d(x,y)=\inf\{1/n\colon y\in B_n(x)\}.$$
To argue that the topology of $X$ is generated by the distance function $d$,
take any point $x\in X$ and any open set $U\subset X$ with $x\in U$.
There is an $n\in\N$ such that $B_n(x)\subset U$.
By Lemma \ref{beobachtung}, $B_d(x,1/n)\subset B_n(x)\subset U$.

If $X$ is first countable, then the sets $B_n(x)$ may be assumed to be open.
To argue that $d$ is well-behaved, fix $x\in X$ and $n\in\mathbb N$.
By Lemma \ref{beobachtung}, $B_{n+1}(x)\subset B_d(x,1/n).$
In particular, $x\in  Int(B_d(x,1/n))$.
\end{proof}

Recall that a topological space $X$ is
\begin{itemize}
\item {\em sequential}
if for each non-closed subset $A\subset X$ there is a sequence $\{a_n\}_{n\in\w}\subset A$ that converges to a point $a\in\bar A\setminus A$;
\item {\em Fr\'echet-Urysohn} if for any subset $A\subset X$ and a point $a\in\bar A\setminus A$ there is a sequence $\{a_n\}_{n\in\w}\subset A$ that converges to $a$.
\end{itemize}
It is rather easy to see from these definitions that
a topological space $X$ is Fr\'echet-Urysohn if and only if each subspace of $X$ is sequential.
Distance spaces, being weakly first countable, are sequential, see \cite{SD}.
It is tempting to think that each subspace of a distance space is itself a distance space,
and that in effect all distance spaces are Fr\'echet-Urysohn, but this reasoning is wrong.

\begin{example}
There is a Hausdorff symmetric (non-well-behaved) distance space $(X,d)$
that is not Fr\'echet-Urysohn and which contains
a subspace $Z\subset X$ that is not sequential,
whose subspace topology is not generated
by the restricted distance function $d|(Z\times Z)$.
\end{example}
\begin{proof}
Let $K=\{1/n\colon n\in\N\}$, $A=K\times K$, $Z=\{(0,0)\}\cup A$, $X=Z\cup(K\times\{0\})$.
Let $d\colon X\times X\to[0,1]$ be the Euclidean metric
except $$d\big((0,0),(1/n,1/m)\big)=d\big((1/n,1/m),(0,0)\big)=1.$$
Notice that $(X,d)$ is not well-behaved because $Int(B(a,0.9))$ is empty.
Notice that $a\in\clo A\setminus A$, but no sequence of points in $A$ converges to $a$.
This means that $X$ is not Fr\'echet-Urysohn
and that $Z=\{a\}\cup A$ is not sequential.
Although $\{a\}$ is not a relatively open subset of $Z$,
we have $B_d(a,0.9)\cap Z=\{a\}$, which shows that the topology of $Z$
is not generated by $d|(Z\times Z)$.
\end{proof}

\begin{theorem}\label{hausdorffdistance}
For a Hausdorff distance space,
being Fr\'echet-Urysohn is equivalent to being first countable.
\end{theorem}
\begin{proof}
Recall that all first countable topological spaces are Fr\'echet-Urysohn.

Let $(X,d)$ be a Fr\'echet-Urysohn Hausdorff distance space.
By Theorem \ref{summary2}, the compression map $\pi\colon\Cob(X,d)\to X$
is a quotient continuous surjection.
According to \cite[2.4.F(c)]{En}, a quotient surjection onto
a Fr\'echet-Urysohn Hausdorff space is hereditarily quotient.
So $\pi$ is hereditarily quotient.
By Theorem \ref{summary2}, $(X,d)$ is well-behaved and thus first countable.
\end{proof}

\section{The functors $\Cob$ and $\Cob^\omega$}\label{Eco}

Recalling Definition \ref{cobomega}, let us introduce the following notation.
\begin{definition}{\em
Let $(X,d)$ be a distance space and
let $\pi_0^1\colon\Cob(X,d)\to X$ denote the compression map.
Let $\pi^\w_1\colon\Cob^\w(X,d)\to\Cob^1(X,d)$ denote the limit projection,
given by $\pi^\w_1(x_1,\ldots,x_n,\ldots)=x_1\in\Cob(X,d)$.

Let $\pi^\omega_0=\pi^1_0\circ\pi^\w_1$ be called the limit compression map
$\pi^\w_0\colon\Cob^\w(X,d)\to X$.
}\end{definition}

\begin{theorem}\label{limitprojection}
For any distance space $(X,d)$ the limit compression map\\
$\pi^\w_0\colon\Cob^\w(X,d)\to X$ is a continuous monotone and quotient surjection.\\
Moreover, it is hereditarily quotient iff $(X,d)$ is well-behaved.
\end{theorem}
\begin{proof}
The limit compression map $\pi^\omega_0=\pi^1_0\circ\pi^\w_1$
is continuous as the composition of two continuous functions.
It is surjective because $\pi^\omega_0((x)_{n=1}^\infty)=x$ for all $x\in X$.

By Theorem~\ref{summary2}, the maps
$\pi^{n+1}_n\colon\Cob^{n+1}(X,d)\to\Cob^n(X,d)$ are continuous monotone hereditarily quotient surjections.
Therefore, by \cite[Theorem 9]{Puzio}, the limit compression map
$\pi^\w_1\colon\Cob^1(X,d)\to\Cob^1(X,d)$ is hereditarily quotient.
By the Corollary to \cite[Theorem 11]{Puzio}, it is also monotone.

Since $\pi^\w_1$ is hereditarily quotient and $\pi^\w_1,\pi^1_0$ are both monotone,
by Lemma \ref{stillmonotone}, $\pi^\omega_0=\pi^1_0\circ\pi^\w_1$ is monotone.

By Theorem \ref{summary2}, $\pi^1_0$ is quotient;
moreover, it is hereditarily quotient if $(X,d)$ is well-behaved.
Now, $\pi^\omega_0$ is the composition of quotient maps so it is quotient.
Moreover, if $(X,d)$ is well-behaved, 
$\pi^\omega_0$ is the composition of hereditarily quotient maps,
so it is hereditarily quotient, by Lemma \ref{herquocompo}.

Assume now that $\pi^\omega_0$ is hereditarily quotient.
To show that $(X,d)$ is well-behaved, take any $x\in X$ and any $r\in(0,1/2)$.
Note that $x\in\Cob^n(X,d)$ for all $n\in\N$.
Let $$U_1=\Cob^\omega(X,d)\cap\Big(B_{\rho_1}(x,r)\times\prod_{n=2}^\infty\Cob^n(X,d)\Big),$$
$$U_2=\{z\in\Cob^\omega(X,d)\colon\pi^\omega_0(z)=x\}\setminus\{(x)_{n=1}^\infty\}.$$
Notice that both of these sets are open and that $(\pi^\omega_0)^{-1}(x)\subset U=U_1\cup U_2$.
Since $\pi^\omega_0$ is hereditarily quotient, $x\in Int(\pi^\omega_0(U))$.
For any $(z_n)_{n=1}^\infty\in U_1$, we have
$z_1\in B_{\rho_1}(x,r)$ and consequently
$\pi^\omega_0(z)=\pi^1_0(\pi^\omega_1(z))=\pi^1_0(z_1)\in B_d(x,r).$
In effect, $\pi^\omega_0(U_1)\subset B_d(x,r)$.
This means that $x\in Int(\pi^\omega_0(U))\subset Int(B_d(x,r))$.
\end{proof}

We are now ready for our strongest result.

\begin{theorem}\label{strongest}
Each (connected) first countable space is the image of a (connected) complete economical metric space
under a continuous monotone hereditarily quotient map.
Each (connected) weakly first countable space is the image of a (connected) complete economical metric space
under a continuous monotone quotient map.
\end{theorem}
\begin{proof}
Let $X$ be a first countable space.
By Theorem \ref{1stcountablewell}, the topology of $X$ is generated by
a well-behaved distance function $d\colon X\times X\to[0,1]$.
By Theorem \ref{limitprojection},
$\pi_0^\w\colon\Cob^\w(X,d)\to X$ is a continuous monotone hereditarly quotient surjection.
By Theorem \ref{complete_economical},
$(\Cob^\w(X,d),\rho_\infty)$ is a complete economical metric space.
By Theorem \ref{summary10},
$\Cob^\w(X,d)$ is connected if $X$ is connected.
We argue analogously for weakly first countable spaces.
\end{proof}

\section{Applications of the cobweb functor}\label{applications}

In this section we shall apply the cobweb construction
to obtain a non-constant continuous locally extremal function
defined on a connected complete metric space.
We define a function $f\colon X\to\IR$ to be {\em locally extremal}
if each point $x\in X$ is a point of local maximum or local minimum of $f$.
In \cite{Ser} Sierpi\'nski proved that each continuous locally extremal function $f\colon\IR\to \IR$ is constant.
This result was generalized in \cite{BGN} to locally extremal functions $f\colon X\to\IR$
defined on connected metrizable spaces $X$ of density $\dens(X)<\mathfrak c$.
In this situation it is natural to ask if the result of Sierpi\'nski holds true for connected metrizable spaces of arbitrary density.
A counterexample to this problem was constructed by Le Donne and Fedeli in \cite{DF2},
and independently by the authors in \cite{BVW}.
Here we present another example based on the cobweb construction.

\begin{theorem}\label{extremal}
There is a connected complete metric space $Z$
with cardinality and density both equal to $\mathfrak c$
and a monotone continuous hereditarily quotient surjection $f\colon Z\to(0,1)$
that has a local extremum at every point.
\end{theorem}
\begin{proof}
Let $X=(0,1)\times\{0,1\}$.
Let $d\colon X\times X\to[0,1]$ be a symmetric distance function on $X$
such that for all $t\in(0,1)$ and $\e\in(0,t)\cap(0,1-t)$,
\begin{enumerate}[(i)]
\item $d\big((t,0),(t,1)\big)=0$
\item $d\big((t,0),(t-\e,1)\big)=\e$
\item $d\big((t,1),(t+\e,0)\big)=\e$
\item $d(a,b)=1$ otherwise.
\end{enumerate}

Let $Z=\Cob(X,d)$ with the natural cobweb metric $\rho$
induced from the complete graph over $X$.
Then, by Theorem \ref{summary},
$Z$ is a complete metric space such that $|Z|=\dens(Z)=\mathfrak c$.

Let $p\colon X\to(0,1)$ be given by $p(t,y)=t$.
Let $\pi\colon\Cob(X,d)\to X$ denote the compression map.
Let $f=p\circ\pi$. Then $f\colon Z\to(0,1)$ is a surjection
that is locally constant at each $z\in Z\setminus X$.

Thanks to (i) the fibers of $f$ are connected, so $f$ is monotone.

Notice that for all $t\in(0,1)$ and $\e\in(0,t)\cap(0,1-t)$,
$$p(\pi(B_\rho((t,0),\e)))=p(B_d((t,0),\e))=(t-\e,t],\text{ by (ii) and (iv),}$$
$$p(\pi(B_\rho((t,1),\e)))=p(B_d((t,1),\e))=[t,t+\e),\text{ by (iii) and (iv)}.$$
This means, in particular, that $f$ has a local maximum
at each point $(t,0)\in Z$ and a local minimum at each point $(t,1)\in Z$.
So $f$ has a local extremum at every point.
Moreover, it follows that $f$ is continuous.

Take any $t\in(0,1)$ and any open set $U\subset Z$ such that $f^{-1}(t)\subset U$.
Then $(t,0),(t,1)\in U$.
Since $U$ is open, there is an $\e\in(0,\frac{1}{2})$ such that
$A=B_\rho((t,0),\e)\cup B_\rho((t,1),\e)\subset U$.
But then $(t-\e,t+\e)\cap(0,1)=f(A)\subset f(U)$,
so $t\in Int(f(U))$.
By Definition \ref{defs}, $f$ is hereditarily quotient.
By Lemma \ref{ZconnectedifXconnected}, $Z$ is connected.
\end{proof}

Recall that open surjections and closed surjections are always quotient maps.
However, a continuous mononotone hereditarily quotient surjection
may be arbitrarily irregular as shown in the following theorem,
which makes use of the cobweb construction.

\begin{theorem}[Krzysztof Omiljanowski]\label{Omiljanowski}
Let $(X,d)$ be a (connected) metric space and let $E\subset X$ be an arbitrary subset.
There is a (connected) complete metric space $(Z,\rho)$ and a monotone continuous
hereditarily quotient surjection $f\colon Z\to X$ such that
for every $x\in X$ there are points $a,b\in Z$ such that $f(a)=f(b)=x$ and
for all $\e\in(0,\frac{1}{2})$,
$$f\big(B_\rho(a,\e)\big)=\{x\}\cup\big(B_d(x,\e)\cap E\big)$$
$$f\big(B_\rho(b,\e)\big)=\{x\}\cup\big(B_d(x,\e)\setminus E\big).$$
\end{theorem}
\begin{proof}
Let $Y=X\times\{0,1\}$.
Let $r\colon Y\times Y\to[0,\infty)$ be a nonsymmetric distance function on $Y$
defined by
$$r\big((x,a),(y,b)\big)=\begin{cases}
0&\text{if }x=y,\\
d(x,y)&\text{if }a=0\wedge b=0\wedge y\in E,\\
d(x,y)&\text{if }a=1\wedge b=0\wedge y\not\in E,\\
1&\text{otherwise.}\end{cases}$$

Let $Z=\Cob(Y,r)$ with the natural cobweb metric induced from
the complete graph over $Y$. Thus $Z$ is a complete metric space.

Let $p\colon Y\to X$ be given by $p(x,a)=x$.
Let $\pi\colon\Cob(Y,r)\to Y$ denote the compression map.
Let $f=p\circ\pi$.
Then $f\colon Z\to X$ is a surjection that is locally constant
at each $Z\setminus Y$, and thus continuous at these points.

Since $r((x,0),(x,1))=0$, the fibers of $f$ are connected, so $f$ is monotone.

Notice that for all $x\in X$ and $\e\in(0,\frac{1}{2})$,
$$\pi(B_\rho((x,0),\e))=B_r((x,0),\e)=\{(x,0)\}\cup\{(y,0)\colon y\in E\wedge d(x,y)<\e\},$$
$$f(B_\rho((x,0),\e))=\{x\}\cup(B_d(x,\e)\cap E),$$
$$\pi(B_\rho((x,1),\e))=B_r((x,1),\e)=\{(x,1)\}\cup\{(y,0)\colon y\not\in E\wedge d(x,y)<\e\},$$
$$f(B_\rho((x,1),\e))=\{x\}\cup(B_d(x,\e)\setminus E).$$
This means, in particular, that $f$ is continuous at each $z\in Y\subset Z$.

Moreover, $f$ is hereditarily quotient by the same argument as in the proof of Theorem \ref{extremal}.
By Lemma \ref{ZconnectedifXconnected}, $Z$ is connected if $X$ is connected.
\end{proof}

\section{Appendix: Quotient maps}\label{appendix}

Quotient maps play a fundamental role in our crucial results.
In this section we collect a number of definitions and known facts.

\begin{definition}\label{defs}
Let $X$ and $Y$ be topological spaces.
Then a function $f\colon X\to Y$ is
\begin{enumerate}[(1)]
\item {\em monotone} if
$f^{-1}(y)$ is connected for each $y\in Y$;
\item {\em quotient} if $A$ is open in $Y$ whenever $f^{-1}(A)$ is open in $X$;
\item {\em hereditarily quotient}
if $f^{-1}(y)\subset U\implies y\in Int(f(U))$
\\for every open set $U\subset X$ and every point $y\in Y$.
\end{enumerate}
\end{definition}

Usually, a function $f\colon X\to Y$ is defined to be hereditarily quotient
if for every $Z\subset Y$ the restriction $f|f^{-1}(Z)\colon f^{-1}(Z)\to Z$ is quotient,
which is easily implied by property Definition \ref{defs}(4), as shown in Theorem \ref{eng121}.
However, according to \cite[2.4.F(a)]{En}, these two conditions are equivalent
for a continuous surjection.

When dealing with continuous surjections in our theorems,
each time we show a function to be hereditarily quotient
we actually show the stronger property and each time
we assume a function to be hereditarily quotient,
we actually need the stronger property.
So we wrote the definition like that
to avoid the necessity of going through the proof of equivalence.

\begin{theorem}\label{eng121}
Let $X,Y$ be topological spaces and let $f\colon X\to Y$.
Suppose that for every $y\in Y$ and every open set $U\subset X$ we have
$f^{-1}(y)\subset U\implies y\in Int(f(U))$.
Then for every $Z\subset Y$ the restriction $f|f^{-1}(Z)\colon f^{-1}(Z)\to Z$ is quotient.
\end{theorem}
\begin{proof}
Take any $A\subset Z\subset Y$
such that $f^{-1}(A)$ is open in $f^{-1}(Z)$.
There is an open set $U\subset X$ such that $f^{-1}(A)=f^{-1}(Z)\cap U$.
Notice that $Z\cap f(U)\subset A$.

Since $f^{-1}(A)\subset U$, for every $y\in A$, $f^{-1}(y)\subset U$.
Therefore $A\subset Int(f(U))$, by assumption.
Now, $A\subset Int(f(U))\cap Z\subset f(U)\cap Z\subset A$.
So $A=Int(f(U))\cap Z$ is open in $Z$.
We showed that $f|f^{-1}(Z)$ is quotient for every $Z\subset Y$.
\end{proof}


The following lemma is used in Section \ref{cobweb}.

\begin{lemma}
\label{ZconnectedifXconnected}
Suppose that $f\colon X\to Y$ is monotone and quotient.
Then X is connected if f(X) is connected.
\end{lemma}
\begin{proof}
Let $U$ be a clopen subset of $X$.
Since $f$ is monotone, $U=f^{-1}(f(U))$.
Since $f$ is quotient, $f(U)$ is open in $Y$.
Since $X\setminus U$ is also clopen,
it follows by analogy that $f(X\setminus U)$ is open, too.
Now, the sets $f(U)$, $f(X\setminus U)$ are disjoint,
because the clopen sets $U$, $X\setminus U$ contain whole fibers.
The connected space $f(X)=f(U)\cup f(X\setminus U)$
is a union of two disjoint open subsets.
Thus $U=X$ or $U=\emptyset$, showing that $X$ is connected.
\end{proof}

The following lemmas prepare for Puzio's Theorem \ref{Puzio},
which is used in Section \ref{iterated}.

\begin{lemma}\label{supermonotone}
If $f\colon X\to Y$ is monotone and hereditarily quotient
then $f^{-1}(E)$ is connected whenever $E\subset Y$ is connected.
\end{lemma}
\begin{proof}
The restriction $f|{f^{-1}(E)}\colon f^{-1}(E)\to E$ is quotient
and evidently monotone.
By Lemma~\ref{ZconnectedifXconnected},
$f^{-1}(E)$ is connected if $E$ is connected.
\end{proof}

\begin{lemma}\label{stillmonotone}
If $f\colon X\to Y$ and $g\colon Y\to Z$ are monotone
and $f$ is hereditarily quotient,
then $g\circ f$ is monotone.
\end{lemma}
\begin{proof}
Since $g^{-1}(z)$ is connected
and $f$ is monotone and hereditarily quotient,
by Lemma~\ref{supermonotone},
$(g\circ f)^{-1}(z)=f^{-1}(g^{-1}(z))$ is connected.
\end{proof}

\begin{lemma}\label{herquocompo}
If $f\colon X\to Y$ and $g\colon Y\to Z$ are hereditarily quotient
then their composition $g\circ f\colon X\to Z$ is hereditarily quotient.
\end{lemma}
\begin{proof}
Take any open set $U\subset X$ and any point $z\in Z$
such that $(g\circ f)^{-1}(z)\subset U$. Then
$(g\circ f)^{-1}(z)=f^{-1}(g^{-1}(z)))=\bigcup\{f^{-1}(y)\colon y\in g^{-1}(z)\}\ \subset\ U.$
\\Since $f$ is hereditarily quotient, we have $y\in Int(f(U))$
for each $y\in g^{-1}(z)$.
\\Thus $g^{-1}(z)\subset Int(f(U))$.
Now, since $g$ is hereditarily quotient,
\\$z\in Int(g(Int(f(U))))\subset Int((g\circ f)(U)).$
\end{proof}

\begin{theorem}[E. Puzio, 1972, \cite{Puzio}]\label{Puzio}
Let $X_1,X_2,\ldots$ be a sequence of connected spaces.
Suppose that each function $f_n\colon X_{n+1}\to X_n$
is a continuous monotone hereditarily quotient surjection.
Then the inverse limit of this system
$$X=\big\{(x_n)_{n=1}^\infty\in\prod_{n=1}^\infty X_n\colon
(\forall n\in\mathbb N)\big(x_n=f_n(x_{n+1})\big)\big\}$$
is connected.
\end{theorem}
\begin{proof}
Theorem 11 in \cite{Puzio}
requires the additional assumption that the functions
$$f_n\circ f_{n+1}\circ\ldots\circ f_m\colon X_{m+1}\to X_n$$
are monotone hereditarily quotient surjections
for all $n<m\in\mathbb N$.
But this follows from our assumptions
thanks to Lemma~\ref{stillmonotone} and Lemma~\ref{herquocompo}.
So $X$ is connected.
\end{proof}

The little example $[0,1)\cup\{2\}\rightarrow[0,1]$
shows that the assumption that $f$ is hereditarily quotient
is needed in all these lemmas.

\section{Problems and comments}

Let us make some basic observations about the notion of an economical metric.

To see that the Cantor set is economically metrizable,
consider the Cantor cube $2^\N=\{0,1\}^\N$ endowed with the metric
$d(x,y)=\max\{|x(n)-y(n)|/n\colon n\in\N\}$,
which assumes only countably many values.
On the other hand,
the Euclidean metric $d(x,y)=|x-y|$ on the Cantor ternary set $C$
is not economical, because
$\dens(C)=\aleph_0<\mathfrak c=|C|=|d(\{0\}\times C)|\le|d(C\times C)|$.
Every metrizable space containing a copy of the Cantor set
admits a metric that is not economical.
Indeed, recall that given a metrizable space $X$ and a closed subset $M\subset X$,
every admissable metric on $M$ can be extended to an admissable metric on $X$,
see \cite[4.5.21(c)]{En}.
So if $X$ contains a copy of the Cantor set, $M\subset X$,
we may choose a metric for $M$ that is not economical
and extend it to the whole $X$.

Let us state without proof that a metric space $(X,d)$
is economical if it satisfies at least one of these conditions:
(1) $d$ assumes only countably many values, $|d(X\times X)|\le\aleph_0$,
(2) $(X,d)$ is a scattered space (contains no nonempty dense-in-itself subsets),
(3) $d$ satisfies the strong triangle inequality (= $d$ is an ultrametric),
$d(x,z)\le\max\{d(x,y),d(y,z)\}\text{ for all }x,y,z\in X,$.
\\
\\
Let us notice that the axiom of choice was not used in proving that
for every connected metric space $(X,d)$ there is a connected
complete metric space whose separablewise components
form a quotient space homeomorphic to $X$.
However, to prove that the complete economical metric space $\Cob^\omega(X,d)$
is actually connected, we made use of Puzio's Theorem 11 in \cite{Puzio}
which uses the axiom of choice. This leaves Problem 1 of \cite{MW2} still open.
\begin{problem}
Is there a nonseparably connected complete metric space in ZF?
\end{problem}

In \cite{MW2}, Morayne and W\'ojcik constructed a nonseparably connected metric group,
which is an example of a homogeneous nonseparably connected metric space.

\begin{problem}
Is there a nonseparably connected complete metric space that is homogeneous?
\end{problem}

\begin{problem}
Is there a locally connected nonseparably connected metric space?
\end{problem}

We would like to thank Pawe{\l} Krupski for his topological seminar
at which we had a chance to present our work and improve it
greatly along the lines suggested by Krzysztof Omiljanowski.

We rely on \cite{Puzio} in our Theorem \ref{summary10}
to argue that the inverse limit is connected,
and in Theorem \ref{limitprojection} to argue that the limit projection
is hereditarily quotient.
In Section \ref{distance} on distance spaces
we quote some classical results from Engelking's book \cite{En}
and refer to \cite{SD} to argue that distance spaces are sequential.
Besides that, our proofs do not require any other sources.

\end{document}